\documentclass[final,twoside,11pt]{entics}
\usepackage{enticsmacro}
\usepackage{graphicx}
\usepackage[all]{xy}

\usepackage{cite}
\usepackage{amsmath}
\allowdisplaybreaks

\def\conf{ISDT 2022} 	
\volume{2}			
\def\volu{2}			
\def\papno{6}			
\def\mydoi{{\normalshape  \href{https://doi.org/10.46298/entics.10339}{10.46298/entics.10339}}}
\def\copynames{M. Liu, Y. Li} 
\def\runtitle{Quantaloidal completions of rrder-enriched categories and their applications}

\def\lastname{Liu and Li}
\begin{document}
\begin{frontmatter}
  \title{Quantaloidal Completions of Order-enriched Categories and Their Applications}
  \author{Min Liu
  \thanksref{ALL}\thanksref{myemail}}   
   \author{Yulin Li
   \thanksref{coemail}}        
  \thanks[ALL]{This work is
supported by the National Natural Science Foundation of China (Grant
No. 11871320), Natural Science Basic Research Program of Shaanxi
(Program No. 2022JM-032) and the Fundamental Research Funds for the
Central Universities (Grant No. 300102122109).}   
   \thanks[myemail]{Email: \href{mailto:liumin@chd.edu.cn} {\texttt{\normalshape
        liumin@chd.edu.cn}}}
  \address
  {School of Sciences\\Chang'an University\\
     Xi'an, P R China}
  \thanks[coemail]{Email:  \href{mailto:1942348167@qq.com} {\texttt{\normalshape
        1942348167@qq.com}}}
\begin{abstract}
By introducing the concept of quantaloidal completions for an
order-enriched category, relationships between the category of
quantaloids and the category of order-enriched categories are
studied. It is proved that quantaloidal completions for an
order-enriched category can be fully characterized as compatible
quotients of the power-set completion. As applications, we show that
a special type of injective hull of an order-enriched category is
the MacNeille completion; the free quantaloid over an order-enriched
category is the Down-set completion.
\end{abstract}
\end{frontmatter}
\section{Introduction}\label{intro}
An order-enriched category is a locally small category such that the
hom-sets are partially ordered sets and composition of morphisms
preserve order in both  variables. An order-enriched category with
only one object can be viewed as a partially ordered semigroup. Thus
order-enriched categories can be viewed as categorical
generalization of partially ordered semigroups.  Several works
devoted to this subject are from computer science
\cite{MartinHoareHe,Wand1979}, especially with strong background of
the study of programming languages. In 1979, M. Wand studied
fixed-point constructions in order-enriched categories, which
extended Scott's result based on continuous lattices. Note that an
order-enriched category in the sense of \cite{Wand1979} means a
category with hom-sets not only ordered but also with certain
completeness.  Later, M. Smyth and G. Plotkin considered solving
recursive domain equations in this framework \cite{SmythPlotkin}. In
2007, this ideal was further extended to the framework of
bicategories\cite{CattaniFiore}. In 1991, C. E. Martin, C. A. R.
Hoare and He Jifeng studied pre-adjunctions in order
enriched-categories \cite{MartinHoareHe}. In \cite{MartinHoareHe},
the concepts of lax functors, natural transformations and
pre-adjunctions are studied with the purpose to explain their
understanding of programming languages. We also note that an
order-enriched category in the sense of \cite{MartinHoareHe} means a
category with hom-sets preordered. These works are all devoted to
study special kind of order-enriched categories. There are little
works devoted to study on them systematically.

A quantaloid $\mathcal{Q}$
\cite{AbramskyVickers,Rosenthal1991,Rosenthal1995,Rosenthal1996} is
a category enriched in the symmetric monoidal closed category {\bf
Sup} of complete lattices and morphisms that preserve arbitrary
sups. Just as every complete lattice is a special partially ordered
set, every quantaloid is a special order-enriched categories. A
quantaloid with only one object is a quantale \cite{Rosenthal1990},
thus quantaloids are naturally viewed as quantales with many
objects. Quantaloids were studied by Pitts \cite{Pitts} in
investigating distributive categories of relations and topos theory
under the name of sup-lattice enriched categories. In
\cite{AbramskyVickers} quantaloids are studied in order to include a
notion of type on the processes. Quantaloids and their applications
were further developed in the monograph \cite{Rosenthal1996}. In
recent years, Quantaloid-enriched categories received considerable
attention
\cite{CleHof2009,Garner,laishen,LiuZhao20131,Hohle2014,Shen2016,ShenTao2016,ShenTholen2016,Stubbe20051,Stubbe20052,Stubbe20053,Stubbe2007,Stubbe2010,Stubbe2014}.

 The process of completion is a classic approach
to study ordered structures. Various completion methods for ordered
structures are developed with different characteristics
\cite{MacNeille,Banaschewski,BanaschewskiNelson,Erne,HanZhao,Rasouli,XiaZhaoHan}.
Relationships between order-enriched categories and quantaloids have
not received enough attention, though they have similar backgrounds
and close relations. Inspired by research on completion methods for
ordered semigroups and their
applications\cite{HanZhao,Lambek2012,Rosenthal1991,XiaZhaoHan}, this
paper is devoted to study quantaloidal completions of order-enriched
categories and their applications.

The contents of the paper are arranged as follows. Section 2 lists
some preliminary notions and results about order-enriched categories
and quantaloids. In Section 3, based on compatible nuclei on
quantaloids, quantaloidal completions for an order-enriched category
are fully characterized as compatible quotients of the power-set
completion. In Section 4, two aspects of applications of
quantaloidal completions are given. It is proved that the injective
hull of an order-enriched category with respect to a special kind of
morphisms is the MacNeille completion; the free quantaloid over an
order-enriched category is the Down-set completion.


\section{Preliminaries on order-enriched categories and quantaloids}
For category theory, we refer to \cite{AdamekHS,Kelly1982}. Let
$\mathcal{C}_{0}$ be the class of objects of a category
$\mathcal{C}$. $\mathcal{C}(a, b)$ denotes the hom-set for $a, b\in
\mathcal{C}_{0}$. For $a\in \mathcal{C}_{0}$, $
1_{a}$ denotes the
identity on $a.$

\begin{definition} (\cite{MartinHoareHe})\ \  An \emph{order-enriched category} is a locally small category $\mathcal{A}$
such that:
\begin{enumerate}
\item[(1)] for $a, b\in \mathcal{A}_{0}$, the hom-set $\mathcal{A}(a, b)$
is a poset,

\item[(2)] composition of morphisms of $\mathcal{A}$ preserves order in
both variables. 
\end{enumerate}
\end{definition}

\begin{definition}  (\cite{ZhangPaseka2019})\ \ Let
$\mathcal{C}$, $\mathcal{D}$ be order-enriched categories. A
\emph{lax semifunctor} $F:\mathcal{C}\rightarrow \mathcal{D}$ is
given by functions $F: \mathcal{C}_{0}\rightarrow \mathcal{D}_{0}$
and $F_{a,b}: \mathcal{C}(a, b)\rightarrow \mathcal{D}(Fa, Fb)$ for
all $a, b\in \mathcal{C}_{0}$ such that $F_{a,b}$ is
order-preserving and $(Fg)\circ (Ff)\leq F(g\circ f)$ for all $a, b,
c\in \mathcal{C}_{0}$, $f\in \mathcal{C}(a, b), g\in \mathcal{C}(b,
c)$. A \emph{lax functor}  $F:\mathcal{C}\rightarrow \mathcal{D}$ is
a lax semifunctor such that $1_{Fa}\leq F(1_{a})$ for all $a\in
\mathcal{C}_{0}$. A \emph{2-functor} $F: \mathcal{C}\rightarrow
\mathcal{D}$ is a functor such that
$$F_{a,b}: \mathcal{C}(a, b)\rightarrow \mathcal{D}(Fa, Fb)$$
is order-preserving for all $a, b\in \mathcal{C}_{0}.$
\end{definition}

A \emph{quantaloid} $\mathcal{Q}$ \cite{Rosenthal1996} is a category
enriched in the symmetric monoidal closed category {\bf Sup} of
complete lattices and morphisms that preserve arbitrary sups.
 In elementary terms:

\begin{definition}  (\cite{Rosenthal1991})\ \  A \emph{quantaloid}
is a locally small category $\mathcal{Q}$ such that:
\begin{enumerate}
\item[(1)] for $a, b\in \mathcal{Q}_{0}$, the hom-set $\mathcal{Q}(a, b)$
is a complete lattice,

\item[(2)] composition of morphisms of $\mathcal{Q}$ presevers sups in both
variables.
\end{enumerate}
\end{definition}

 In this paper,
$\mathcal{Q}$ always denotes a small quantaloid, and
$\mathcal{Q}_{0}$ denotes the set of its objects. The identity
$\mathcal{Q}$-arrow on $q\in \mathcal{Q}_{0}$ will be denoted by
$1_{q}.$ The greatest element of the complete lattice
$\mathcal{Q}(p, q)$ will be denoted by $\top_{p,q}$. For a
$\mathcal{Q}$-arrow $u: p\longrightarrow q$, we denote the domain
and the codomain of $u$ by $\mathrm{dom} (u)$ and $\mathrm{cod}(u)$,
respectively. Given $\mathcal{Q}$-arrows $u: p\longrightarrow q,$
$v: q\longrightarrow r,$ the corresponding adjoints induced by the
compositions $-\circ u: \mathcal{Q}(q, r)\longrightarrow
\mathcal{Q}(p, r)$ and $v\circ -: \mathcal{Q}(p, q)\longrightarrow
\mathcal{Q}(p, r)$ are denoted by $u\rightarrow_{l}-$ and
$v\rightarrow_{r}-$ respectively.

For more details on quantaloids, we refer to
\cite{Rosenthal1991,Rosenthal1996}.

\begin{definition} (\cite{Rosenthal1991})\ \ Let $\mathcal{Q}$,
$\mathcal{S}$ be quantaloids. A \emph{quantaloidal homomorphism} $F:
\mathcal{Q}\rightarrow \mathcal{S}$ is a functor such that
$$F: \mathcal{Q}(X, Y)\rightarrow \mathcal{S}(FX, FY)$$
is sup-preserving for all $X, Y \in \mathcal{Q}_{0}.$
\end{definition}

A \emph{quantaloidal isomorphism} is a quantaloidal homomorphism
such that it is bijictive on objects and hom-sets.

\begin{example}  Let $\mathcal{A}$ be an order-enriched
category.
\begin{enumerate}
\item[(1)]  $\mathcal{P}(\mathcal{A})$ is a quantaloid~\cite{Rosenthal1991}. The objects of $\mathcal{P}(\mathcal{A})$
are those of $\mathcal{A}$. For $a, b\in \mathcal{A}$, the hom-set
$\mathcal{P}(\mathcal{A})(a, b)$=${\mathcal{P}}(\mathcal{A}(a, b))$,
the power set of the hom-set $\mathcal{A}(a, b)$. For $S\in
\mathcal{P}(\mathcal{A})(a, b), T\in \mathcal{P}(\mathcal{A})(b,
c)$, $T\circ S=\{g\circ f\mid g\in T, f\in S\}.$

\item[(2)]  $\mathcal{D}(\mathcal{A})$ is a quantaloid. The objects of
$\mathcal{D}(\mathcal{A})$ are those of $\mathcal{A}$. For $a, b\in
\mathcal{A}$, the hom-set $\mathcal{D}(\mathcal{A})(a,
b)$=${\mathcal{D}}(\mathcal{A}(a, b))$, the set of down
sets\footnote{A set $D$ in a poset $P$ is a down set, if
$D=\downarrow D$, where $\downarrow D=\{x\mid \exists d\in D,\
\mathrm{s.\ t.}\ x\leq d\}$.} of the hom-set $\mathcal{A}(a, b)$.
For $S\in \mathcal{D}(\mathcal{A})(a, b), T\in
\mathcal{D}(\mathcal{A})(b, c)$, $T\circ S=\downarrow\{g\circ f\mid
g\in T, f\in S\}.$ We note that $\downarrow 1_{a}\in
{\mathcal{D}}(\mathcal{A}(a, a))$ is the identity morphism.
\end{enumerate}
\end{example}

\begin{definition} (\cite{Rosenthal1991})  Let $\mathcal{Q}$ be
a quantaloid. A \emph{quantaloidal nucleus} is a lax functor
$j:\mathcal{Q}\rightarrow\mathcal{Q}$, which is the identity on the
objects of $\mathcal{Q}$ and such that the maps $j_{a, b}:
\mathcal{Q}(a, b)\rightarrow \mathcal{Q}(a, b)$ satisfy:
\begin{enumerate}
\item[(1)] $f\leq j_{a,b}(f)$ for all $f\in \mathcal{Q}(a, b)$,

\item[(2)] $j_{a,b}(j_{a,b}(f))=j_{a,b}(f)$ for all $f\in \mathcal{Q}(a,
b)$,

\item[(3)] $j_{b,c}(g)\circ j_{a,b}(f)\leq j_{a,c}(g\circ f)$ for all $g\in
\mathcal{Q}(b, c)$ $f\in \mathcal{Q}(a, b)$.
\end{enumerate}
\end{definition}

For a quantaloidal nucleus  $j$  on a quantaloid $\mathcal{Q}$, let
$\mathcal{Q}_{j}$ be the bicategory with the same objects as
$\mathcal{Q}$ and $\mathcal{Q}_{j}(a, b)=\{f\in\mathcal{Q}(a, b)\mid
j_{a,b}(f)=f\}$ for $a,  b\in (\mathcal{Q}_{j})_{0}$.  Composition
in $\mathcal{Q}_{j}$ is defined as follows: $g\circ_{j}f=j_{a,
c}(g\circ f)$ for $f\in \mathcal{Q}_{j}(a, b)$, $g\in
\mathcal{Q}_{j}(b, c)$.

\begin{proposition}(\cite{Rosenthal1991}) If $j$ is a
quantaloidal nucleus on a quantaloid $\mathcal{Q}$, then
$\mathcal{Q}_{j}$ is a quantaloid and $j:\mathcal{Q}\rightarrow
\mathcal{Q}_{j}$ is a quantaloidal homomorphism.
\end{proposition}

\begin{proposition} (\cite{Rosenthal1991}) Let $\mathcal{S}$ be
a subcategory of  a quantaloid $\mathcal{Q}$, which contains all the
objects of $\mathcal{Q}$. Then, $\mathcal{S}$ is a quotient
quantaloid of the form $\mathcal{Q}_{j}$ for some quantaloidal
nucleus $j$ iff
\begin{enumerate}
\item[(1)] each hom-set $\mathcal{S}(a, b)$ is closed under infs, and

\item[(2)] if $f\in \mathcal{S}(a, c)$, then $g\rightarrow_{l}f\in
\mathcal{S}(b, c)$ for all $g\in \mathcal{Q}(a, b)$ and
$h\rightarrow_{r}g\in \mathcal{S}(a ,b)$ for all $h\in
\mathcal{Q}(b, c).$ 
\end{enumerate}
\end{proposition}

\section{Quantaloidal completions of order-enriched categories}
In order to study quantaloidal completions of order-enriched
categories, let us begin with the concept of a compatible nucleus on
a quantaloid.

\begin{definition}Let $\mathcal{A}$ be an order-enriched
category, $j:\mathcal{P}(\mathcal{A})\rightarrow
\mathcal{P}(\mathcal{A})$ a quantaloidal nucleus. $j$ is said to be
\emph{compatible} if for $a, b\in \mathcal{A}_{0}$,
$f\in\mathcal{A}(a, b)$, we have $j_{a, b}(\{f\})=\downarrow f$.
\end{definition}

\begin{definition} Let $\mathcal{A}$ be an order-enriched
category, $\mathcal{Q}$ a quantaloid, $F: \mathcal{A}\rightarrow
\mathcal{Q}$ a 2-functor. The pair $(F, \mathcal{Q})$ is said to be
a \emph{quantaloidal completion} of $\mathcal{A}$, if the following
conditions are satisfied:
\begin{enumerate}
\item[(1)] $F: \mathcal{A}_{0}\rightarrow \mathcal{Q}_{0}$ is bijective,

\item[(2)] $F_{a,b}: \mathcal{A}(a, b)\rightarrow \mathcal{Q}(Fa, Fb)$ is
an order embedding for all $a ,b\in \mathcal{A}_{0}$,

\item[(3)] for every $a ,b\in \mathcal{A}_{0}$ and $f\in \mathcal{Q}(Fa,
Fb)$, there exists $U_{f}\subseteq \mathcal{A}(a, b)$ such that
$f=\bigvee F(U_{f})$. 
\end{enumerate}
\end{definition}

\begin{theorem}If $j$ is a compatible nucleus on an
order-enriched category $\mathcal{A}$, then $(F_{j},
\mathcal{P}(\mathcal{A})_{j})$ is a quantaloidal completion of
$\mathcal{A}$, where $F_{j}: \mathcal{A}\rightarrow
\mathcal{P}(\mathcal{A})_{j}$ is defined as follows:
\begin{enumerate}
\item[(1)] $F_{j}: \mathcal{A}_{0}\rightarrow
(\mathcal{P}(\mathcal{A})_{j})_{0}$ is the identity map,

\item[(2)] $F_{j}(f)=\ \downarrow\! f$ for every $f\in\mathcal{A}(a, b)$, $a
,b\in \mathcal{A}_{0}.$ 
\end{enumerate}
\end{theorem}

\begin{proof} By definition, $F_{j}: \mathcal{A}_{0}\rightarrow
(\mathcal{P}(\mathcal{A})_{j})_{0}$ is bijective, and $F_{j}:
\mathcal{A}(a, b)\rightarrow \mathcal{P}(\mathcal{A})_{j}(a, b)$ is
an order-embedding. For $S\in \mathcal{P}(\mathcal{A})_{j}(a, b)$,
we have $S=j_{a, b}(S)=j_{a, b}(\bigcup_{f\in S}
\{f\})=\bigvee_{f\in S}^{\mathcal{P}(\mathcal{A})_{j}(a, b)}j_{a,
b}(\{f\})=\bigvee_{f\in S}^{\mathcal{P}(\mathcal{A})_{j}(a,
b)}\downarrow f=\bigvee_{f\in S}^{\mathcal{P}(\mathcal{A})_{j}(a,
b)}F_{j}(f).$ This completes the proof. \end{proof}

Corresponding to several classical completion methods of posets and
ordered semigroups, we can obtain a series of  compatible nucleus.
We leave detail to the reader.

\begin{example} (Down-set completion) Let $\mathcal{A}$ be an
order-enriched category. Define a lax functor $\downarrow:
\mathcal{P}(\mathcal{A})\rightarrow \mathcal{P}(\mathcal{A})$ as
follows:
\begin{enumerate}
\item[(1)] $\downarrow: \mathcal{P}(\mathcal{A})_{0}\rightarrow
\mathcal{P}(\mathcal{A})_{0}$ is the identity map,

\item[(2)] $\downarrow_{a, b}(S)=\downarrow S$ for $S\in
\mathcal{P}(\mathcal{A})(a, b)$, $a, b\in
\mathcal{P}(\mathcal{A})_{0}.$
\end{enumerate}
Then $\downarrow$ is a compatible nucleus. The quotient
corresponding to $\downarrow$ is $\mathcal{D}(\mathcal{A})$.
\end{example}

\begin{example}(MacNeille completion) Let $\mathcal{A}$ be an
order-enriched category. Define a lax functor $\mathrm{cl}:
\mathcal{P}(\mathcal{A})\rightarrow \mathcal{P}(\mathcal{A})$ as
follows:

(1) $\mathrm{cl}: \mathcal{P}(\mathcal{A})_{0}\rightarrow
\mathcal{P}(\mathcal{A})_{0}$ is the identity map,

(2) $\mathrm{cl}_{a, b}(S)=\{f\in \mathcal{P}(\mathcal{A})(a, b)\mid
\forall g\in \mathcal{P}(\mathcal{A})(a^{\prime}, a), h\in
\mathcal{P}(\mathcal{A})(b, b^{\prime}), k\in
\mathcal{P}(\mathcal{A})(a, b), h\circ S \circ g\subseteq \downarrow
k\ \mathrm{implies}\ h\circ f\circ g\leq k\}$ for $S\in
\mathcal{P}(\mathcal{A})(a, b)$, $a, b\in
\mathcal{P}(\mathcal{A})_{0}.$
\\
Then $\mathrm{cl}$ is a compatible nucleus. \end{example}

\begin{example}(Equivariant completion) Let $\mathcal{A}$ be an
order-enriched category. Suppose $S\subseteq
\mathcal{P}(\mathcal{A})(a, b)$. If the join of $S$ exists and is
preserved by composition, i.e., $f\circ (\bigvee S)=\bigvee (f\circ
S)$, $(\bigvee S)\circ g=\bigvee (S\circ g)$ whenever the
composition is well-defined, then $\bigvee S$ is said to be an
\emph{equivariant join} with respect to $S$. Clearly, every $f\in
\mathcal{P}(\mathcal{A})(a, b)$ is an equivariant join respect to
$\downarrow f.$ If $k$ is an equivariant join with respect to $S$,
then $g\circ k$ (resp., $k\circ h$) is an equivariant join with
respect to $g\circ S$ (resp., $S\circ h$), whenever the composition
is well-defined. For  $S\subseteq \mathcal{P}(\mathcal{A})(a, b)$,
let
$$S^{EJ}=\{f\in \mathcal{P}(\mathcal{A})(a, b)\mid \exists T\subseteq S, \mathrm{s. t.}\ f=\bigvee T\ \mathrm{is\ an\ equivariant\ join\ with\ respect\ to\ } T\}.$$
Let $EJ(\mathcal{A})$ be the subcategory of $\mathcal{A}$, which
contains all the objects of $\mathcal{A}$. The hom-sets
$$EJ(\mathcal{A})(a, b)=\{S\in \mathcal{D}(\mathcal{A})(a,
b)\mid S=S^{EJ} \}.$$ Then $EJ(\mathcal{A})$ is a quotient of
$\mathcal{A}$ such that $\downarrow f\in EJ(\mathcal{A})(a, b)$ for
every $f\in \mathcal{P}(\mathcal{A})(a, b)$. Consequently, the
corresponding quantaloidal nucleus is compatible. \end{example}

For an order-enriched category $\mathcal{A}$, $CN(\mathcal{A})$
denotes the class of all compatible nuclei on
$\mathcal{P}(\mathcal{A})$, $QC(\mathcal{A})$ denotes the set of all
quantaloidal completions of $\mathcal{A}$.

Let $\mathcal{A}$ be an order-enriched category, $(F,
\mathcal{Q})\in QC(\mathcal{A})$. Define $j_{(F, \mathcal{Q})}:
\mathcal{P}(\mathcal{A})\rightarrow \mathcal{P}(\mathcal{A})$ as
follows:
\begin{enumerate}
\item[(1)] $j_{(F, \mathcal{Q})}: \mathcal{P}(\mathcal{A})_{0}\rightarrow
\mathcal{P}(\mathcal{A})_{0}$ is the identity map,

\item[(2)] $j_{(F, \mathcal{Q})}(S)=\{f\in \mathcal{A}(a, b)\mid  F(f)\leq
\bigvee_{g\in S}F(g)\}$ for every $S\in \mathcal{P}(\mathcal{A})(a,
b), a, b\in \mathcal{A}_{0}$.
\end{enumerate}

\begin{lemma}Let $\mathcal{A}$ be an order-enriched category,
$(F, \mathcal{Q})\in QC(\mathcal{A})$. Then $j_{(F, \mathcal{Q})}$
is a compatible nucleus on $\mathcal{P}(\mathcal{A})$. \end{lemma}

\begin{proof} By definition, $j_{(F, \mathcal{Q})}:
\mathcal{P}(\mathcal{A})_{0}\rightarrow
\mathcal{P}(\mathcal{A})_{0}$ is bijective, $j_{(F,
\mathcal{Q})}:\mathcal{P}(\mathcal{A})(a, b)\rightarrow
\mathcal{P}(\mathcal{A})(a, b)$ is order preserving and increasing
for all $a, b\in \mathcal{A}_{0}$. Suppose  $S\in
\mathcal{P}(\mathcal{A})(a, b), f\in j_{(F, \mathcal{Q})}(j_{(F,
\mathcal{Q})}(S))$. Then, $F(f)\leq \bigvee_{g\in j_{(F,
\mathcal{Q})}(S)}F(g)$. For every $g\in j_{(F, \mathcal{Q})}(S),$ we
have $F(g)\leq\bigvee_{k\in S}F(k)$. Thus, $F(f)\leq \bigvee_{k\in
S}F(k)$. Consequently, $f\in j_{(F, \mathcal{Q})}(S)$.
 So we can conclude that $j_{(F, \mathcal{Q})}\circ j_{(F, \mathcal{Q})}=j_{(F,
 \mathcal{Q})}$. Thus, $j_{(F,
\mathcal{Q})}:\mathcal{P}(\mathcal{A})(a, b)\rightarrow
\mathcal{P}(\mathcal{A})(a, b)$ is a closure operator for every $a,
b\in \mathcal{A}_{0}$.

Suppose $K\in \mathcal{P}(\mathcal{A})(b, c)$, $S\in
\mathcal{P}(\mathcal{A})(a, b)$. Then $j_{(F, \mathcal{Q})}(K)\circ
j_{(F, \mathcal{Q})}(S)=\{g\circ f\mid  g\in \mathcal{A}(b, c), f\in
\mathcal{A}(a, b), F(g)\leq\bigvee_{k\in K}F(k),
F(f)\leq\bigvee_{t\in S}F(t)\}$. If $F(g)\leq\bigvee_{k\in K}F(k),
F(f)\leq\bigvee_{t\in S}F(t)$, then $F(g\circ f)\leq \bigvee_{k\in
K, t\in S}F(k)\circ F(t)=\bigvee_{k\in K, t\in S}F(k\circ
t)\leq\bigvee_{p\in K\circ S}F(p).$ Thus, $j_{(F,
\mathcal{Q})}(K)\circ j_{(F, \mathcal{Q})}(S)\subseteq j_{(F,
\mathcal{Q})}(K\circ S)$.

 For $f_{0}\in \mathcal{A}(a, b)$, by the fact that $F: \mathcal{A}(a, b)\rightarrow \mathcal{Q}(F(a),
 F(b))$ is an order embedding, we have $j_{(F,
\mathcal{Q})}(\{f_{0}\})=\{f\in \mathcal{A}(a, b)\mid  F(f)\leq
F(f_{0})\}=\downarrow  f_{0}$.

So we can conclude that $j_{(F, \mathcal{Q})}$ is a compatible
nucleus on $\mathcal{P}(\mathcal{A})$. \end{proof}

\begin{theorem} Let $\mathcal{A}$ be an order-enriched category,
$(F, \mathcal{Q})\in QC(\mathcal{A})$. Then $\mathcal{Q}$ is
quantaloidal isomorphism to $\mathcal{P}(\mathcal{A})_{j_{(F,
\mathcal{Q})}}$. \end{theorem}

\begin{proof} Let $F^{-1}: \mathcal{Q}_{0}\rightarrow
\mathcal{A}_{0}$ be the inverse of the map $F:
\mathcal{A}_{0}\rightarrow \mathcal{Q}_{0}$. Define $G:
\mathcal{Q}\rightarrow \mathcal{P}(\mathcal{A})_{j_{(F,
\mathcal{Q})}}$ as follows:

(1) $G(a)=F^{-1}(a)$ for every $a\in \mathcal{Q}_{0}$,

(2) $G(p)=\{f\in\mathcal{A}(F^{-1}(c), F^{-1}(d))\mid  F(f)\leq p\}$
for every $p\in \mathcal{Q}(c, d)$.
\\
Then $G: \mathcal{Q}_{0}\rightarrow
(\mathcal{P}(\mathcal{A})_{j_{(F, \mathcal{Q})}})_{0}$ is bijective.
For $f\in j_{(F, \mathcal{Q})}(G(p))$, we have
$F(f)\leq\bigvee_{g\in G(p)}F(g)\leq p$, thus $f\in G(p)$. Thus,
$j_{(F, \mathcal{Q})}(G(p))\subseteq G(p)$. Consequently,
$G(p)=j_{(F, \mathcal{Q})}(G(p))\in \mathcal{P}(\mathcal{A})_{j_{(F,
\mathcal{Q})}}$. Thus, $G$ is well-defined.

Suppose $a, b\in \mathcal{Q}_{0}, S\subseteq \mathcal{Q}(a, b)$.
Then $G(\bigvee S)=\{f\in \mathcal{A}(F^{-1}(a), F^{-1}(b))\mid
F(f)\leq \bigvee S\}$,
\begin{align*}
\bigvee_{t\in S}^{\mathcal{P}(\mathcal{A})_{j_{(F,
\mathcal{Q})}}}G(t)&=j_{(F, \mathcal{Q})}\left(\bigcup_{t\in S}G(t)\right)\\
&=j_{(F, \mathcal{Q})}\left(\bigcup_{t\in S}\{g\in
\mathcal{A}(F^{-1}(a), F^{-1}(b))\mid
F(g)\leq t\}\right)\\
&=j_{(F, \mathcal{Q})}\{g\in \mathcal{A}(F^{-1}(a), F^{-1}(b))\mid
\exists t\in S, \mathrm{\ s. t.\ }
F(g)\leq t\}\\
&=\{f\in \mathcal{A}(F^{-1}(a), F^{-1}(b))\mid  F(f)\leq
\bigvee\{F(g)\mid  g\in \mathcal{A}(F^{-1}(a), F^{-1}(b)), \exists
t\in S,\\ &\mathrm{\ \ \ \ \ \  s. t.\ } F(g)\leq t\}\}.
\end{align*}
For $s_{0}\in S,$ we have
\begin{align*}
s_{0}&=\bigvee\{F(g)\mid g\in \mathcal{A}(F^{-1}(a), F^{-1}(b)),
F(g)\leq s_{0}\}\\
&\leq \bigvee\{F(g)\mid g\in \mathcal{A}(F^{-1}(a), F^{-1}(b)),
\exists t\in S, \mathrm{\ s. t.\ } F(g)\leq t\}.
\end{align*}
Thus, $\bigvee S\leq \bigvee\{F(g)\mid g\in \mathcal{A}(F^{-1}(a),
F^{-1}(b)), \exists t\in S, \mathrm{\ s. t.\ } F(g)\leq t\}$, whence
$G(\bigvee S)\leq\bigvee_{t\in S}^{\mathcal{P}(\mathcal{A})_{j_{(F,
\mathcal{Q})}}}G(t)$. The inverse inequality holds trivially.
Therefore, $G(\bigvee S)=\bigvee_{t\in
S}^{\mathcal{P}(\mathcal{A})_{j_{(F, \mathcal{Q})}}}G(t)$.

For $a\in \mathcal{Q}_{0}$, we have $G(1_{a})=\{f\in
\mathcal{A}(F^{-1}(a), F^{-1}(a))\mid  F(f)\leq 1_{a}\}=\downarrow
1_{G(a)}$, which is the identity in
$\mathcal{P}(\mathcal{A})_{j_{(F, \mathcal{Q})}}.$

Suppose $f\in \mathcal{Q}(a, b), g\in \mathcal{Q}(b, c)$. Then
\begin{align*}
G(g)\circ_{j_{(F, \mathcal{Q})}}G(f)&=j_{(F, \mathcal{Q})}(G(g)\circ
G(f))\\
&=\{t\in \mathcal{A}(F^{-1}(a), F^{-1}(c))\mid  F(t)\leq
\bigvee\{F(h)\mid h\in G(g)\circ G(f)\}\}.
\end{align*}
Since,
\begin{align*}
&\bigvee\{F(h)\mid h\in G(g)\circ G(f)\}\\
=&\bigvee\{F(t_{2}\circ t_{1})\mid  t_{1}\in \mathcal{A}(F^{-1}(a),
F^{-1}(b)), t_{2}\in \mathcal{A}(F^{-1}(b), F^{-1}(c)), F(t_{1})\leq
f, F(t_{2})\leq g\}\\
=&\left(\bigvee\{F(t_{2})\mid  t_{2}\in \mathcal{A}(F^{-1}(b),
F^{-1}(c)), F(t_{2})\leq g\}\right)\circ\left(\bigvee\{F(t_{1})\mid
t_{1}\in \mathcal{A}(F^{-1}(a), F^{-1}(b)),  F(t_{1})\leq
f\}\right)\\
=&g\circ f,
\end{align*}
we have $G(g)\circ_{j_{(F, \mathcal{Q})}}G(f)=\{t\in
\mathcal{A}(F^{-1}(a), F^{-1}(c))\mid F(t)\leq g\circ f\}=G(g\circ
f)$.

So we can conclude that $G$ is a quantaloidal homomorphism.

Suppose $p_{1}, p_{2}\in \mathcal{Q}(c, d)$ with
$G(p_{1})=G(p_{2})$. Then $p_{1}=\bigvee F(G(p_{1}))=\bigvee
F(G(p_{2}))=p_{2}$. Thus, $G: \mathcal{Q}(c, d)\rightarrow
(\mathcal{P}(\mathcal{A})_{j_{(F, \mathcal{Q})}})(F^{-1}(c),
F^{-1}(d))$ is injective for all $c, d\in \mathcal{Q}_{0}$.

Suppose $S\in \mathcal{P}(\mathcal{A})_{j_{(F, \mathcal{Q})}}(a,
b)$. Then $S\subseteq \mathcal{A}(a, b)$. For every
$f\in\mathcal{A}(a, b)$, we have
$G(F(f))=\{g\in\mathcal{A}(F^{-1}(a), F^{-1}(b))\mid F(g)\leq
F(f)\}=\{g\in\mathcal{A}(F^{-1}(a), F^{-1}(b))\mid g\leq
f\}=\downarrow f$. Thus, $S=j_{(F, \mathcal{Q})}(S)=j_{(F,
\mathcal{Q})}\left(\bigcup_{f\in S}\{f\}\right)=\bigvee_{f\in
S}^{\mathcal{P}(\mathcal{A})_{j_{(F, \mathcal{Q})}}(a, b)}j_{(F,
\mathcal{Q})}(\{f\})=\bigvee_{f\in
S}^{\mathcal{P}(\mathcal{A})_{j_{(F, \mathcal{Q})}}(a, b)}\downarrow
f=\bigvee_{f\in S}^{\mathcal{P}(\mathcal{A})_{j_{(F,
\mathcal{Q})}}(a, b)} G(F(f))=G\left(\bigvee_{f\in
S}^{\mathcal{Q}(F(a), F(b))}F(f)\right)$. Thus, $G: \mathcal{Q}(a,
b)\rightarrow (\mathcal{P}(\mathcal{A})_{j_{(F,
\mathcal{Q})}})(F^{-1}(a), F^{-1}(b))$ is surjective. Therefore, $G$
is a quantaloidal isomorphism. \end{proof}

As a combination of the above results, we obtain that quantaloidal
completions of an order-enriched category $\mathcal{A}$ are
completely determined by compatible quantaloidal nuclei on
$\mathcal{P}(\mathcal{A})$.

\begin{theorem} Let $\mathcal{A}$ be an order-enriched category.
Then $(F, \mathcal{Q})$ is a quantaloidal completion of
$\mathcal{A}$ if and only if there is a compatible nucleus $j$ on
$\mathcal{P}(\mathcal{A})$ such that $\mathcal{Q}$ is quantaloidal
isomorphism to $\mathcal{P}(\mathcal{A})_{j}$.\end{theorem}

\section{Applications}
In this section, we shall give two kinds of applications for the
quantaloidal completions of order-enriched categories.

\subsection{Injective constructs of
order-enriched categories}

 Let
$\mathbf{O}$-$\mathbf{Cat}_{l}$ be the category of order-enriched
categories and lax semifunctors. Let $\mathcal{E}^{ls}_{\leq}$ be
the class of all lax semifunctors in $\mathbf{O}$-$\mathbf{Cat}_{l}$
satisfying the following conditions:
\begin{enumerate}
\item[(1)] $F: \mathcal{C}_{0}\rightarrow \mathcal{D}_{0}$ is bijective;

\item[(2)] $F(f_{1})\circ F(f_{2})\circ\cdots\circ F(f_{n})\leq F(f)$ implies
$f_{1}\circ f_{2}\circ\cdots\circ f_{n}\leq f$ for $f_{1}\circ f_{2}\circ\cdots\circ
f_{n}, f\in C(a, b)$, $a, b\in \mathcal{C}_{0}$.
\end{enumerate}

\begin{lemma} In the category $\mathbf{O}$-$\mathbf{Cat}_{l}$,
every retract of a quantaloid is a quantaloid. \end{lemma}

\begin{proof} Let $\mathcal{S}$ be a retract of a quantaloid
$\mathcal{Q}$. Then there exist lax semifunctors $I:
\mathcal{S}\rightarrow \mathcal{Q}$ and $F: \mathcal{Q}\rightarrow
\mathcal{S}$ such that $F\circ I=\mathrm{id}_{\mathcal{S}}$. Suppose
$S, T\in \mathcal{S}_{0}$. Then $\mathcal{S}(X, T)$ is a retract of
$\mathcal{Q}(IX, IY)$. By the fact that $\mathcal{Q}(IX, IY)$ is a
complete lattice, we can deduce that $\mathcal{S}(X, Y)$ is a
complete lattice and $F(\bigvee I(A))$ is the least upper bound of
$A$ in $\mathcal{S}(X, Y)$. Suppose $A\subseteq \mathcal{S}(X, Y)$,
$g\in \mathcal{S}(Y, Y^{\prime})$, $t\in \mathcal{S}(X^{\prime},
X)$. Then $g\circ (\bigvee A)$ is an upper bound of $g\circ A$. If
$h$ is an upper bound of $g\circ A$, then $I(g)\circ \bigvee_{f\in
A}I(f)=\bigvee_{f\in A}(I(g)\circ I(f))\leq \bigvee_{f\in
A}(I(g\circ f))\leq I(h).$ Thus, $h=FI(h)\geq F\left(I(g)\circ
\bigvee_{f\in A}I(f)\right)\geq FI(g)\circ F\left(\bigvee_{f\in
A}I(f)\right)=g\circ(\bigvee A)$. Thus, $g\circ (\bigvee
A)=\bigvee(g\circ A)$. Similarly, we have $(\bigvee A)\circ
t=\bigvee (A\circ t)$. Therefore, $\mathcal{S}$ is a quantaloid.
\end{proof}

\begin{theorem} Let $\mathcal{A}$ be an order-enriched
category. Then $\mathcal{A}$ is $\mathcal{E}^{ls}_{\leq}$-injective
in $\mathbf{O}$-$\mathbf{Cat}_{l}$  if and only if $\mathcal{A}$ is
a quantaloid. \end{theorem}

\begin{proof} Suppose $\mathcal{Q}$ is a quantaloid, $H:
\mathcal{S}\rightarrow\mathcal{T}$ a morphism in
$\mathcal{E}^{ls}_{\leq}$, and  $F:
\mathcal{S}\rightarrow\mathcal{Q}$  a morphism in
$\mathbf{O}$-$\mathbf{Cat}_{l}$. Define $G: \mathcal{T}\rightarrow
\mathcal{Q}$ as follows:

(1) $GX=FH^{-1}(X)$, $\forall X\in \mathcal{T}_{0}$;

(2) $G(g)=\bigvee\{F(f_{1})\circ F(f_{2})\circ \cdots \circ
F(f_{n})\mid  H(f_{1})\circ H(f_{2})\circ \cdots \circ  H(f_{n})\leq
g, f_{1}\circ f_{2}\circ \cdots \circ  f_{n}\in \mathcal{S}(H^{-1}X,
H^{-1}Y)\}$ for $g\in \mathcal{T}(X, Y)$, $X, Y\in \mathcal{T}_{0}$.
\\
Then $G: \mathcal{T}(X, Y)\rightarrow \mathcal{Q}(GX, GY)$ is
order-preserving for $X, Y\in \mathcal{T}_{0}$. Suppose $g_{1}\in
\mathcal{T}(X, Y)$, $g_{2}\in \mathcal{T}(Y, Z)$. Since composition
in a quantaloid distribute over arbitrary joins, we can deduce that
$G(g_{2})\circ G(g_{1})\leq G(g_{2}\circ g_{1})$. Thus $G:
\mathcal{T}\rightarrow \mathcal{Q}$ is a lax semifunctor. For $X\in
\mathcal{S}$, we have $GH(X)=FH^{-1}H(X)=F(X)$. For $f\in
\mathcal{S}(X, Y)$, we have $H(f)\in \mathcal{T}(HX, HY)$. By the
fact that $F$ is  a morphism in $\mathbf{O}$-$\mathbf{Cat}_{l}$, we
can deduce that $GH(f)=\bigvee\{F(f_{1})\circ F(f_{2})\circ \cdots
\circ F(f_{n})\mid  H(f_{1})\circ H(f_{2})\circ \cdots \circ
H(f_{n})\leq H(f), f_{1}\circ f_{2}\circ \cdots \circ  f_{n}\in
\mathcal{S}(X, Y)\}=F(f)$. Thus, $G: \mathcal{T}\rightarrow
\mathcal{Q}$ is a morphism in   $\mathbf{O}$-$\mathbf{Cat}_{l}$ such
that $GH=F$. So we can conclude that $\mathcal{Q}$ is
$\mathcal{E}^{ls}_{\leq}$-injective.

Conversely, suppose $\mathcal{A}$ is
$\mathcal{E}^{ls}_{\leq}$-injective in
$\mathbf{O}$-$\mathbf{Cat}_{l}$.  Define $F: \mathcal{A}\rightarrow
\mathcal{D}(\mathcal{A})$ as follows:

(1) $F: \mathcal{A}_{0}\rightarrow \mathcal{D}(\mathcal{A})_{0}$ is
the identity map;

(2) $F(f)=\downarrow f$ for $f\in \mathcal{A}(a, b)$, $a, b\in
\mathcal{A}_{0}$.
\\
Then its routine to check that $F\in \mathcal{E}^{ls}_{\leq}$. Thus,
for the identity functor $\mathrm{id}_{\mathcal{A}}:
\mathcal{A}\rightarrow \mathcal{A}$, there is a lax semifunctor $G:
\mathcal{D}(\mathcal{A})\rightarrow \mathcal{A}$ such that
$GF=\mathrm{id}_{\mathcal{A}}$. So, $\mathcal{A}$ is a quantaloid,
as it is a retract of the quantaloid $\mathcal{D}(\mathcal{A})$.
\end{proof}

Let $\mathcal{A}$ be an order-enriched category. Define $\eta:
\mathcal{A}\rightarrow \mathcal{P}(\mathcal{A})_{\mathrm{cl}}$ as
follows:
\begin{enumerate}
\item[(1)] $\eta: \mathcal{A}_{0}\rightarrow
(\mathcal{P}(\mathcal{A})_{\mathrm{cl}})_{0}$ is the identity map;

\item[(2)] $\eta(f)=\ \downarrow\! f$ for $f\in \mathcal{A}(a, b)$, $a, b\in
\mathcal{A}_{0}$.
\end{enumerate}
Then it is routine to check that $\eta$ is a lax semifunctor and it
is $\mathcal{E}^{ls}_{\leq}$-essential. As the proof is quite
similar to that of Theorem 5.8 in \cite{Lambek2012}, we leave it to
the reader.

\begin{theorem} Let $\mathcal{A}$ be an order-enriched
category. Then $\mathcal{P}(\mathcal{A})_{\mathrm{cl}}$ is an
$\mathcal{E}^{ls}_{\leq}$-injective hull of $\mathcal{A}$ in
$\mathbf{O}$-$\mathbf{Cat}_{l}$. \end{theorem}

\subsection{Free quantaloids}

Let $\mathbf{LocSm}$  be the category of locally small categories
and functors between them. Let $\mathbf{O}$-$\mathbf{Cat}$ be the
category of order-enriched categories and 2-functors. Let
$\mathbf{Qtlds}$ be the category of quantaloids and quantaloidal
homomorphisms.

\begin{theorem} The functor $\mathcal{D}:$
$\mathbf{O}$-$\mathbf{Cat}$ $\rightarrow  \mathbf{Qtlds}$ is left
adjoint to the forgetful functor $\mathbf{Qtlds} \rightarrow$
$\mathbf{O}$-$\mathbf{Cat}$. \end{theorem}

\begin{proof} Let $\mathcal{A}$ be an order-enriched category. Define
$\eta: \mathcal{A}\rightarrow
  \mathcal{D}(\mathcal{A})$ as follows:

(1) $\eta: \mathcal{A}_{0}\rightarrow (\mathcal{D}(\mathcal{A})
)_{0}$ is the identity map;

(2) $\eta(f)=\downarrow f$ for $f\in \mathcal{A}(a, b)$, $a, b\in
\mathcal{A}_{0}$.
\\
Then $\eta$ is a 2-functor in $\mathbf{O}$-$\mathbf{Cat}$.

Suppose that $\mathcal{Q}$ is a quantaloid and that $F:
\mathcal{A}\rightarrow \mathcal{Q}$ is a 2-functor in
$\mathbf{O}$-$\mathbf{Cat}$. Define $\bar{F}:
\mathcal{D}(\mathcal{A})\rightarrow\mathcal{Q}$ as follows:

(1) $\bar{F}(a)=F(a)$ for every $a\in (\mathcal{D}(\mathcal{A})
)_{0}$;

(2) $\bar{F}(S)=\bigvee\{F(f)\mid f\in S\}$ for every $S\in
\mathcal{D}(\mathcal{A})(a, b)$.
\\
For $a\in (\mathcal{D}(\mathcal{A}))_{0}$, we have
$\bar{F}(\downarrow 1_{a})=\bigvee\{F(f)\mid f\in \downarrow
1_{a}\}=\bigvee\{F(f)\mid f\leq 1_{a}\}=F(1_{a})=1_{F(a)}$. For
$T\in \mathcal{D}(\mathcal{A})(b, c), S\in
\mathcal{D}(\mathcal{A})(a, b)$, we have $\bar{F}(T\circ
S)=\bigvee\{F(h)\mid h\in T\circ S\}=\bigvee\{F(g\circ f)\mid  g\in
T, f\in S\}=\bigvee\{F(g)\mid  g\in T\}\circ \bigvee\{F(f)\mid  f\in
S\}=\bar{F}(T)\circ \bar{F}(S)$. For $S_{i}\in
\mathcal{D}(\mathcal{A})(a, b), i\in I$, we have
$\bar{F}\left(\bigvee_{i\in
I}S_{i}\right)=\bar{F}\left(\bigcup_{i\in
I}S_{i}\right)=\bigvee\{F(f)\mid f\in \bigcup_{i\in
I}S_{i}\}=\bigvee_{i\in I}\bigvee\{F(f)\mid f\in
S_{i}\}=\bigvee_{i\in I}\bar{F}(S_{i})$. Thus, $\bar{F}$ is a
quantaloidal homomorphism. Furthermore, we can check that
$\bar{F}\circ\eta=F$.

Suppose $G: \mathcal{D}(\mathcal{A}\rightarrow\mathcal{Q}$ is a
quantaloidal homomorphism with $G\circ\eta=F$. Then we have

(1) $\forall a\in (\mathcal{D}(\mathcal{A}))_{0}$,
$\bar{F}(a)=\bar{F} (\eta(a))=F(a)=(G\circ \eta)(a)= G(a)$;

(2) $\forall S\in \mathcal{D}(\mathcal{A})(a, b)$,
$\bar{F}(S)=\bar{F}(\bigcup\{\downarrow f\mid f\in
S\})=\bigvee_{f\in S}\bar{F}(\downarrow f)=\bigvee_{f\in
S}\bar{F}(\eta(f))=\bigvee_{f\in S}F(f)=\bigvee_{f\in
S}(G\circ\eta)(f)=\bigvee_{f\in S}G(\downarrow
f)=G\left(\bigvee_{f\in S}\downarrow f\right)=G(S)$. Thus, $\bar{F}:
\mathcal{D}(\mathcal{A}\rightarrow\mathcal{Q}$ is the unique
quantaloidal homomorphism such that $\bar{F}\circ \eta=F$.
\end{proof}

Every locally small category can be viewed as an order-enriched
category with the discrete order on hom-sets. We know
$\mathcal{D}(\mathcal{A})=\mathcal{P}(\mathcal{A})$ for every
locally small category with discrete order on hom-sets. Thus, we can
recover the following results \cite{Rosenthal1991}.

\begin{corollary} The functor $\mathcal{P}:
\mathbf{LocSm}\rightarrow  \mathbf{Qtlds}$ is left adjoint to the
forgetful functor $\mathbf{Qtlds} \rightarrow
\mathbf{LocSm}$.\end{corollary}

\section{Conclusion and some further work}
In this paper, we only considered quantaloidal completions for
order-enriched categories. As order-enriched category with other
completeness have deep applications in domain theory
\cite{MartinHoareHe,SmythPlotkin,Wand1979}, other types of
completions and applications deserve to be developed further.

\begin{ack} The authors would like to thank  the anonymous referees for their valuable comments. \end{ack}

\bibliographystyle{./entics}

\end{document}